\documentclass[twoside,a4paper,10pt]{article}
\usepackage{latexsym,amssymb}
\usepackage{hyperref}
\usepackage{graphicx}
\usepackage{enumerate}
\usepackage[frenchb,english]{babel}
\usepackage[matrix,arrow]{xy}
\usepackage{color}
\def \be{\begin{eqnarray*}}
\def \ee{\end{eqnarray*}}
\def \ben{\begin{enumerate}}
\def \een{\end{enumerate}}
\def \beit{\begin{itemize}}
\def \eeit{\end{itemize}}
\def \bui#1#2{\mathrel{\mathop{\kern 0pt#1}\limits^{#2}}}
\def \buil#1#2{\mathrel{\mathop{\kern 0pt#1}\limits_{#2}}}
\def \bfll{\begin{flushleft}}
\def \efll{\end{flushleft}}
\def \bflr{\begin{flushright}}
\def \eflr{\end{flushright}}

\def \findemo{\hfill$\square$\\}

\def \wit{\widetilde}
\def \wnabla{\wit{\nabla}}
\def \bnabla{\overline{\nabla}}
\def \bg{\overline{g}}

\def \C{\mathbb{C}}

\def \R{\mathbb{R}}

\newtheorem{ethm}{Theorem}[section]

\newtheorem{erem}[ethm]{Note}

\newtheorem{prop}[ethm]{Proposition}

\title{An Obata-type characterization of Calabi metrics on line bundles}
\author{Nicolas Ginoux\footnote{Universit\'e de Lorraine, CNRS, IECL, F-57000 Metz, France, \texttt{nicolas.ginoux@univ-lorraine.fr}
}, Georges Habib\footnote{Lebanese University, Department of Mathematics, P.O. Box 90656 Fanar-Matn, Lebanon, \texttt{ghabib@ul.edu.lb}}, Mihaela Pilca\footnote{Universit\"at Regensburg, Universit\"atstr. 31, 93051 Regensburg, Germany, \texttt{mihaela.pilca@ur.de}} and Uwe Semmelmann\footnote{
Universit{\"a}t Stuttgart,
Pfaffenwaldring 57,
70569 Stuttgart, Germany, \texttt{uwe.semmelmann@mathematik.uni-stuttgart.de}}}
\begin{document}
\maketitle

\begin{center}\begin{tabular}{p{155mm}}
\begin{small}{\bf Abstract.} We characterize those complete K\"ahler manifolds supporting a nonconstant real-valued function with critical points whose Hessian is complex linear, has pointwise two eigenvalues and whose gradient is a Hessian-eigenvector.
\end{small}\\
\end{tabular}\end{center}

\section{Introduction and main results}\label{s:intromainresults}

The celebrated Obata theorem \cite[Theorem A]{Obata62} states that the only complete Riemannian manifold $(M^n,g)$ on which a nonconstant real-valued $C^2$-function $u$ exists whose Hessian satisfies $\nabla^2u=-u\cdot\mathrm{Id}$ is the round sphere up to homothety on the metric.
On K\"ahler manifolds, the Obata theorem was generalized in several ways.
In \cite[Theorem p. 614]{MolzonPinneyMortensen93}, R.~Molzon and K.~Pinney Mortensen proved that a connected complete K\"ahler manifold $(M^{2n},g,J)$ is biholomorphically isometric to $\C\mathrm{P}^n$ with Fubini-Study metric of holomorphic sectional curvature $1$ if and only if there exists a nonconstant function $u\in C^2(M,\R)$ such that
\[\nabla^2 u=-u\cdot\mathrm{Id}+\frac{1}{2}(u-1)\cdot(\mathrm{Id}-\pi)\]
holds on the regular set of $u$, where $\pi$ is the pointwise orthogonal projection onto the subbundle $\mathrm{Span}_\C(\nabla u,J\nabla u)\subset T^\C M$.
Note that in particular $\nabla^2u$ is $J$-invariant and has the pointwise eigenvalues $-u$ and $-\frac{u+1}{2}$.
This result was generalized by A.~Ranjan and G.~Santhanam as follows \cite[Theorem 3]{RanjSant97}: a connected complete K\"ahler manifold $(M^{2n},g,J)$ is biholomorphically isometric to $\C\mathrm{P}^n$ with Fubini-Study metric of holomorphic sectional curvature $1$ if and only if there exists a nonconstant function $u\in C^2(M,\R)$ such that $\nabla^2u$ has pointwise at most two eigenvalues, $-u$ and $-\frac{u+1}{2}$ and that $\nabla u$ is a pointwise eigenvector of $\nabla^2 u$ associated to the eigenvalue $-u$.
Later on, G.~Santhanam considered changing the sign of the eigenvalues of the Hessian of $u$ and prove the following \cite[Theorem 1]{Santhanam07}: 
On a given connected complete K\"ahler manifold $(M^{2n},g,J)$, assume the existence of a nonconstant function $u\in C^2(M,\R)$ with critical points such that $\nabla^2u$ has pointwise two eigenvalues, $u$ and $\frac{u+1}{2}$ and that both $\nabla u$ and $J\nabla u$ are pointwise eigenvectors of $\nabla^2 u$ associated to the eigenvalue $u$.
Then either $u$ has a minimum and in that case $M$ is biholomorphically isometric to $\C\mathrm{H}^n$ with holomorphic sectional curvature $-1$; or $u$ has a maximum and in that case there exists a (real)-$2$-codimensional sub\-ma\-ni\-fold $M_0$ of $M$ such that $M$ is diffeomorphic to $T^\perp M_0$ and where moreover each fibre $T_{x}^\perp M_0$ is isometric to the hyperbolic plane with sectional curvature $-1$.\\

In \cite{GHPSnoncrit}, we weaken the assumptions of the above results even further by leaving the eigenvalues of the Hessian of $u$ arbitrary.
In case the function $u$ has no critical point, we prove the following:

\begin{ethm}[\protect{\cite[Theorem 1.1]{GHPSnoncrit}}]\label{t:cardwp2eigenv}
Let $(\wit{M}^{2n},\wit{g},\wit{J})$ be a connected complete K\"ahler manifold of real dimension $2n\geq4$ carrying a function $u\in C^\infty(\wit{M},\mathbb{R})$ without critical points such that
\begin{enumerate}[$\bullet$]
\item its Hessian $\wnabla^2 u$ is $\wit{J}$-invariant and has pointwise at most two eigenvalues $\lambda$ and $\mu$ and
\item its gradient $\wnabla u$ is a pointwise eigenvector of $\wnabla^2u$ with the eigenvalue $\lambda$.
\end{enumerate}
Then the following claims hold true:
\begin{enumerate}[i)]
\item Case $\mu$ vanishes at one point: then $\mu=0$ and the triple $(\wit{M}^{2n},\wit{g},\wit{J})$ is locally biholomorphically isometric to $(\R_t\times\R_s\times\Sigma,dt^2\oplus\rho(t)^2ds^2\oplus g_{\Sigma})$ for some K\"ahler manifold $(\Sigma,g_{\Sigma})$, for $\rho(t)=|\wnabla u|(t,s,x)$ for all $(t,s,x)\in\R^2\times\Sigma$ and where the K\"ahler structure of $(\R_t\times\R_s\times\Sigma,dt^2\oplus\rho(t)^2ds^2\oplus g_{\Sigma})$ is the one that is naturally induced by that of $(\Sigma,g_{\Sigma})$.
\item Case $\mu$ does not vanish identically:
\begin{enumerate}
\item In case $n>2$: up to changing $u$ into $au+b$ with $a,b\in\R$ and $a\neq0$, the function $u$ must be positive, the K\"ahler manifold $(\wit{M}^{2n},\wit{g},\wit{J})$ is biholomorphically isometric to the doubly-warped product $\left(\R\times M^{2n-1},dt^2\oplus\rho(t)^2\left(\rho'(t)^2\hat{g}_{\hat{\xi}}\oplus\hat{g}_{\hat{\xi}^\perp}\right)\right)$, where $M$ is a level hypersurface of $u$, the triple $(M^{2n-1},\hat{g},\hat{\xi})$ Sasaki and $\rho(t)=\sqrt{u(t,x)}$, for any $(t,x)\in\R\times M$.
\item In case $n=2$: 
up to changing $u$ into $au+b$ with $a,b\in\R$ and $a\neq0$, the function $u$ must be positive, the K\"ahler manifold $(\wit{M}^{2n},\wit{g},\wit{J})$ is biholomorphically isometric to the doubly-warped product $\left(\R\times M^{2n-1},dt^2\oplus\rho(t)^2\left(\rho'(t)^2\hat{g}_{\hat{\xi}}\oplus\hat{g}_{\hat{\xi}^\perp}\right)\right)$, where $M$ is a level hypersurface of $u$, the triple $(M^{2n-1},\hat{g},\hat{\xi})$ is a minimal Riemannian flow that is basic conformally Sasaki and $\rho(t)=\sqrt{u(t,x)}$, for any $(t,x)\in\R\times M$.
\end{enumerate}
\end{enumerate}
Moreover, we have in both cases $\lambda\circ F=\frac{\partial^2 (u\circ F)}{\partial t^2}$ and in the case where $\mu\neq0$ also $\mu=\frac{|\wnabla u|^2}{2u}$.
\end{ethm}

Equivalent formulations of our assumptions as well as the relation to
other classification work are mentioned in the introduction of \cite{GHPSnoncrit}.\\

As noticed in \cite{GHPSnoncrit}, all statements of Theorem \ref{t:cardwp2eigenv} hold locally around every regular point of $u$ on $\wit{M}$ if $(\wit{M}^{2n},\wit{g})$ is not complete or if $u$ has critical points.
The reason is that the identification between $\wit{M}$ and the doubly-warped product is realized by the flow of the normalized gradient vector field of $u$, which exists locally around any regular point of $u$.\\

We look at the case where $u$ still satisfies the same kind of conditions involving its Hessian but is supposed to have critical points.
Note in this case that actually the equation only makes sense on the set of regular points of $u$.
We shall see that, by contrast with \cite{GHPSnoncrit}, the underlying K\"ahler manifold $(\wit{M}^{2n},\wit{g},\wit{J})$ will not always be of warped product form any longer since Calabi metrics on the normal bundle of the critical set of $u$ will turn out to provide further examples:


\begin{ethm}\label{t:cardnormalbundle}
Let $(\wit{M}^{2n},\wit{g},\wit{J})$ be a complete connected K\"ahler manifold.
Assume the existence of a nonconstant function $u\in C^\infty(\wit{M},\R)$ with at least one critical point and with nonnegative $\wit{J}$-invariant Hessian on $\wit{M}$ such that, at any regular point, 
{\begin{enumerate}[$\bullet$]
\item its Hessian $\wnabla^2 u$ has pointwise at most two eigenvalues $\lambda$ and $\mu$ and
\item its gradient $\wnabla u$ is a pointwise eigenvector of $\wnabla^2u$ with the eigenvalue $\lambda$.
\end{enumerate}
Then the critical set $N$ of $u$ is a connected totally geodesic K\"ahler submanifold of $\wit{M}$.
Moreover, the following claims hold true:
\begin{enumerate}[i)]
\item Case $\mu$ vanishes at one regular point of $u$: then $\mu=0$ on $\wit{M}$, the submanifold $N$ has real codimension $2$ and flat normal bundle $T^\perp N\to N$ and the manifold $(\wit{M}^{2n},g,J)$ is biholomorphically isometric to a Calabi construction on the complex line bundle $T^\perp N\to N$.
\item Case $\mu$ does not vanish identically on $\wit{M}$:
\begin{enumerate}
\item If $\mu$ vanishes at a critical point of $u$: then $\mu_{|_N}=0$, the K\"ahler manifold $N$ is Hodge and of real codimension $2$ and the manifold $(\wit{M}^{2n},g,J)$ is biholomorphically isometric to a Calabi construction on the non-flat complex line bundle $T^\perp N\to N$.
\item If $\mu$ does not vanish at any critical point of $u$: then $N=\{p\}$  for some $p\in\wit{M}$ and, up to changing $u$ into $au+b$ with $a,b\in\R$ and $a>0$, the manifold $\wit{M}$ is diffeomorphic to $\mathbb{C}^n$ and $\wit{M}\setminus\{p\}$ with induced K\"ahler structure is isometric to $(0,\infty)\times\mathbb{S}^{2n-1}$ with metric $dt^2\oplus\rho(t)^2(\rho'(t)^2\hat{g}_{\hat{\xi}}\oplus\hat{g}_{\hat{\xi}^\perp})$ where $\rho(t)=\sqrt{u(t,x)}$ for all $(t,x)\in(0,\infty)\times\mathbb{S}^{2n-1}$ and $(\mathbb{S}^{2n-1},\hat{g},\hat{\xi})$ is Sasaki with round metric of constant sectional curvature $1$.
Moreover, $\rho$ must be an odd function, i.e. all even derivatives of $\rho$ must vanish at $0$.\\
Conversely, given any smooth positive odd function $\rho$, the manifold $(0,\infty)\times\mathbb{S}^{2n-1}$ with metric $dt^2\oplus\rho(t)^2(\rho'(t)^2\hat{g}_{\hat{\xi}}\oplus\hat{g}_{\hat{\xi}^\perp})$ can be extended to a complete K\"ahler manifold by adding a point and the function $u:=\rho^2$ satisfies the assumptions above.
\end{enumerate}
\end{enumerate}}
\end{ethm}

We point out that, in case $(\wit{M}^{2n},\wit{g},\wit{J})$ is a Calabi metric on a Hermitian holomorphic line bundle, $u$ as well as $\lambda$ and $\mu$ can be explicitly expressed as functions of the fibrewise length, see Proposition \ref{l:HessianeqonCalabi} below.\\

The article is organized as follows.
In Section \ref{s:prelim}, we recall basics on doubly-warped products and the Calabi construction.
The particular case where $\wit{M}$ is a Calabi construction as well as the proof of Theorem \ref{t:cardnormalbundle} will be carried out in Section \ref{s:proof}.

\section{Preliminaries}\label{s:prelim}

\subsection{K\"ahler doubly-warped products}\label{ss:dwp}

Let us first recall some general facts on product structures induced by smooth functions.
The local version of the following proposition can be found in \cite{Kuehnel}, see in particular \cite[Sec. D]{Kuehnel} or \cite[App.]{GHPSnoncrit}.

\begin{prop}[\protect{\cite[Prop. A.1]{GHPSnoncrit}}]\label{p:warpedprodu}
Let $(M^n,g)$ be a connected complete Riemannian manifold.
Assume that some $u\in C^\infty(M,\R)$ has no critical point on $M^n$ and satisfies $\nabla^2u(\nabla u)=\lambda\nabla u$ for some $\lambda\in C^\infty(M,\R)$, where $\nabla$ denotes the Levi-Civita connection of $(M^n,g)$.
Then the manifold $(M^n,g)$ is isometric to $(\R\times\Sigma,dt^2\oplus g_t)$, where $\Sigma$ is a level hypersurface of $u$ and $(g_t)_{t\in\R}$ is a one-parameter-family of Riemannian metrics on $\Sigma$.
Moreover, the function $u$ only depends on $t$ via that identification.
\end{prop}

The Riemannian manifolds $(M^n,g)$ and $(\R\times\Sigma,dt^2\oplus g_t)$ in Proposition \ref{p:warpedprodu} are identified via the flow of the normalized gradient vector field $\nu:=\frac{\nabla u}{|\nabla u|}$, which in particular preserves the level hypersurfaces of $u$.
If the stronger condition $\nabla^2 u=\lambda\mathrm{Id}$ holds, the metric $g_t$ is of the form $g_t=\rho(t)^2g_{\Sigma}$ for some fixed Riemannian metric $g_{\Sigma}$ on $\Sigma$ and some smooth positive function $\rho$, see e.g. \cite{Kuehnel}.\\

In case $u$ is allowed to have critical points, the same metric splitting holds locally away from the critical points of $u$.
Besides, those cannot be saddle points \cite[Prop. 3]{Santhanam07}.\\

We also recall that a Riemannian doubly-warped product over a triple $(M,\hat{g},\hat{\xi})$ consisting of a Riemannian manifold $(M,\hat{g})$ and a globally defined unit vector field $\hat{\xi}$ is a Riemannian manifold of the form 
\[(\wit{M},\wit{g}):=\left(I\times M,dt^2\oplus\rho(t)^2(\sigma(t)^2\hat{g}_{\hat{\xi}}\oplus\hat{g}_{\hat{\xi}}^\perp)\right),\]
where $\rho,\sigma$ are smooth positive functions on some open interval $I\subset\R$.
In case $(M,\hat{g},\hat{\xi})$ is Sasaki and $\sigma=\rho'$, the doubly-warped product manifold $(\wit{M},\wit{g})$ becomes K\"ahler for the almost complex structure $\wit{J}$ that is naturally induced from the transversal one on $(M,\hat{g},\hat{\xi})$, see \cite[Lemma 2.4]{GHPSnoncrit}.

\subsection{The Calabi construction on a holomorphic Hermitian line bundle}\label{ss:calabiconstr}

We briefly recall the so-called Calabi construction of K\"ahler metrics on Hermitian holomorphic line bundles, which will turn out to provide one of the two main families of examples for our Obata-like equation.
We follow the very clear elementary introduction \cite{Gauduchon13}, see also \cite{ACG01,ChiossiNagy10,HwangSinger02}.\\

Let $L\bui{\longrightarrow}{\pi} N$ be a holomorphic Hermitian line bundle with complex structure $J$ over a K\"ahler manifold $(N,g,J)$ with K\"ahler form $\omega$.
Denote by $\nabla$ the Chern connection on $L$, that is, the unique metric connection on $L$ with $\nabla^{0,1}=\overline{\partial}$, and by $\rho^\nabla\in\Omega^2(N,\R)$ its curvature form, which is defined  for all vector fields $X,Y$ on $N$ by $[\nabla_X,\nabla_Y]-\nabla_{[X,Y]}=:R_{X,Y}^\nabla=i\rho^\nabla(X,Y)\otimes\mathrm{id}_L$.
Let $r\colon L\to\R$ be defined by $r(v):=|v|$ for all $v\in L$, where $|\cdot|$ denotes the Hermitian metric on $L$; the function $r$ is the fibrewise length function, in particular it is smooth on $L^\times:=\{v\in L\,|\,v\neq0\}$ and its vanishing set coincides with the zero section of $L$.
The Calabi ansatz for constructing K\"ahler metrics on the total space of $L$ consists in requiring the real two-form
$$ \overline{\omega}:=h_1\cdot\pi^*\omega+h_2\cdot dr\wedge d^cr,$$
to be the K\"ahler form of some K\"ahler metric $\bg$ on $L$ for suitable positive smooth functions $h_1,h_2$ on $L$ and where $d^cr:=-dr\circ J$.
Note that $J$ preserves $\overline{\omega}$ by definition and holomorphicity of $L$, that is, $\overline{\omega}(JX,JY)=\overline{\omega}(X,Y)$ for all $X,Y\in TL$.
Assuming $h_1$ and $h_2$ to be $\mathbb{U}_1$-invariant, that is, to depend fibrewise only on $r$, the metric $\bg$ is K\"ahler on $L$ if and only if $d\overline{\omega}=0$, condition which is equivalent to $h_1,h_2$ to have vanishing horizontal differential, that is, to depend globally only on $\tau:=r^2$, and to the identities
$$\rho^\nabla=l\omega\qquad\textrm{as well as}\qquad h_1'(\tau)+\frac{l}{2}h_2(\tau)=0$$
for some $l\in\R$.
As a consequence, apart from the situation where $l=0$ and hence $L$ is flat and $h_1$ is constant, the metric $\bg$ is a K\"ahler metric on $L$ only if $(N,g)$ is Hodge, that is, its K\"ahler class is a scalar multiple of an integral class $\alpha$ such that $\frac{1}{m}\alpha\notin H^2(N,\mathbb{Z})$ for every integer $m\geq2$; in that case, $L\to N$ is the Hermitian line bundle associated to some integral multiple of that class $\alpha$ and $h_1$ and $h_2$ are easily obtained solving a linear first-order ODE -- mind that they must be positive though, so that global existence of $\bg$ on $L$ is guaranteed only on a neighbourhood of the zero section.
Note that if $l>0$, then $h_1'<0$ and, unless $\int_0^\infty h_2(\tau)d\tau<\infty$, the metric $\bg$ cannot be complete.\\

From now on $h_1$ and $h_2$ will be considered as functions of $r$, in particular the differential equation they satisfy becomes  $h_1'(r)+lrh_2(r)=0$ for all $r\geq0$.

\section{Proof of the main result}\label{s:proof}

\subsection{Obata-type equation for Calabi metrics}\label{ss:ObataCalabi}

We start with looking for particular functions satisfying our conditions on holomorphic Hermitian line bundles endowed with a metric of Calabi-type.

\begin{prop}\label{l:HessianeqonCalabi}
Let $L\bui{\to}{\pi} N$ be any holomorphic Hermitian line bundle over any K\"ahler manifold $N$.
Let $\overline{\omega}:=h_1(r)\cdot\pi^*\omega+h_2(r)\cdot dr\wedge d^cr$ be the K\"ahler structure on $L$ defined by the Calabi construction associated to $h_1$ and $h_2$, where $\rho^\nabla=l\omega$ and $h_1'(r)+lrh_2(r)=0$ for some $l\in\R$ and where $N$ is assumed to be Hodge when $l\neq0$.
Let $u\colon L\to\R$ be any smooth function depending fibrewise only on $r$ and whose gradient vector field w.r.t. the K\"ahler metric $\bg$ is vertical.
Assume the Hessian of $u$ to have pointwise at most two eigenvalues $\lambda$ and $\mu$ and to have $\mathrm{Span}(\partial_r,J\partial_r)$ as pointwise eigenspace associated to the eigenvalue $\lambda$.\\
Then the function $u$ is either constant or, up to replacing $u$ by $bu+b'$ for some $b,b'\in\R$ with $b\neq0$, given by $u(r)=\int_0^{r}s\cdot h_2(s)ds$ for all $r\geq0$.
Moreover, $\lambda(r)=1+\frac{rh_2'(r)}{2h_2(r)}$ and $\mu(r)=\frac{rh_1'(r)}{2h_1(r)}=-\frac{lr^2 h_2(r)}{4h_1(r)}$ for all $r\geq0$.
\end{prop}

{\it Proof:} 
Recall first the Koszul formula for the Levi-Civita connection $\nabla^M$ on any Riemannian manifold $(M,g)$: for any $X,Y,Z\in\Gamma(TM)$,
\begin{eqnarray}\label{eq:Koszulformula} 
\nonumber g(\nabla_X^M Y,Z)&=&\frac{1}{2}\Big\{X(g(Y,Z))+Y(g(Z,X))-Z(g(X,Y))\\
&&\phantom{\frac{1}{2}\Big\{}+g([X,Y],Z)-g([Y,Z],X)+g([Z,X],Y)\Big\}.
\end{eqnarray}
By definition, the K\"ahler metric $\bg$ associated to $\overline{\omega}$ and $J$ on $L$ is given by
$$\bg=\overline{\omega}(\cdot,J\cdot)=h_1\cdot\pi^* g_N+h_2\cdot(dr\otimes dr+d^cr\otimes d^cr),$$
where we have denoted by $g_N$ the metric on $N$.
As a consequence, the horizontal distribution of the map $\pi$ w.r.t. $\bg$, which is defined by $\mathcal{H}L:=\left(\mathcal{V}L\right)^\perp$ where $\mathcal{V}L_p=T_p\pi^{-1}(\{\pi(p)\})$ is the vertical distribution at any $p\in L$, is given by 
$$\mathcal{H}L=\ker(dr)\cap\ker(d^cr).$$
Note however that $\pi\colon(L,\bg)\to(N,g_N)$ is in general not a Riemannian submersion because of the $h_1$-factor in the metric $\bg$.
It can be shown that actually $\mathcal{H}L$ coincides with the horizontal distribution associated to the Chern connection $\nabla$, that is, for any $p\in L$
$$\mathcal{H}L_p=d_{\pi(p)}s(T_{\pi(p)}N)$$
holds, where $s\in\Gamma(L)$ is any section satisfying $s(\pi(p))=p$ and $(\nabla s)_{\pi(p)}=0$.
As a consequence, for any horizontal lift $X^*$ of some vector field $X$ on $N$, we have $X^*(r)=JX^*(r)=0$.\\
It can also be shown that, for any horizontal vector field $\overline{Z}$ on $L$, that is, $\overline{Z}\in\Gamma(TL)$ with $\overline{Z}_p\in\mathcal{H}L_p$ for all $p\in L$, the Lie bracket $[\overline{Z},\partial_r]$ is again horizontal.
This follows from the flow of $\partial_r$ preserving $\mathcal{H}L$: namely if $(\varphi_t^{\partial_r})_t$ denotes the local flow of $\partial_r$, at least in $L^\times$, then $\varphi_t^{\partial_r}(p)=p+t\partial_r=p+t\frac{p}{|p|}=(1+\frac{t}{r(p)})p$ for all $p\in L^\times$ and $t$ in a neighbourhood of $0\in\R$; since $dr$ vanishes on horizontal vector fields, we have $d_p\varphi_t^{\partial_r}(X)=(1+\frac{t}{r(p)})X$ for every $X\in\mathcal{H}L_p$, in particular $d_p\varphi_t^{\partial_r}(X)\in\mathcal{H}L_p$.
An analogous argument shows that $[J\partial_r,\overline{Z}]$ is horizontal.
Note that, as a particular case, if $\overline{Z}=Z^*$ is the horizontal lift of some vector field $Z$ on $N$, then $[Z^*,\partial_r]$ must also be vertical because of $d\pi([Z^*,\partial_r])=[Z,d\pi(\partial_r)]=0$ and thus $[Z^*,\partial_r]=0$.\\
We compute the different components of the Levi-Civita connection $\bnabla$ of $(L,\bg)$ thanks to the Koszul formula (\ref{eq:Koszulformula}).
Recall that, since we consider $h_1=h_1(r),h_2=h_2(r)$ as functions of $r$, we have $h_1'+rlh_2=0$ because of the change of variable from $\tau$ to $r$.
First,
$$\bg(\bnabla_{\partial_r}\partial_r,\partial_r)=\frac{1}{2}\partial_r(\bg(\partial_r,\partial_r))=\frac{h_2'(r)}{2}.$$
Second, 
\begin{eqnarray*}
\bg(\bnabla_{\partial_r}\partial_r,J\partial_r)&=&\frac{1}{2}\Big(\partial_r(\bg(\partial_r,J\partial_r))+\partial_r(\bg(J\partial_r,\partial_r))-J\partial_r(\bg(\partial_r,\partial_r))+\bg([\partial_r,\partial_r],J\partial_r)\\
&&\phantom{\frac{1}{2}\Big(}-\bg([\partial_r,J\partial_r],\partial_r)+\bg([J\partial_r,\partial_r],\partial_r)\Big)\\
&=&\frac{1}{2}\Big(-\underbrace{J\partial_r(h_2)}_{0}+2h_2dr([J\partial_r,\partial_r])\Big).
\end{eqnarray*}
Since both $J\partial_r$ and $\partial_r$ are tangent to the fibres of $L$, we may for any $p\in L$ consider the flat Levi-Ci\-vi\-ta connection denoted by $\nabla^{\pi(p)}$ of the fibre $(L_{\pi(p)},\langle\cdot\,,\cdot\rangle_{\pi(p)})$ through $p$ with its original Hermitian metric.
By the fact that $\partial_r$ is geodesic w.r.t. $\langle\cdot\,,\cdot\rangle_{\pi(p)}$ and that $\nabla_{\cdot}^{\pi(p)}\partial_r$ is minus the Weingarten endomorphism of distance circles around $0$ (and whose unit normal is $\partial_r$), we have
$$[J\partial_r,\partial_r]_p=\nabla_{J\partial_r}^{\pi(p)}\partial_r-\nabla_{\partial_r}^{\pi(p)}J\partial_r=\frac{J\partial_r}{r(p)}-J(\nabla_{\partial_r}^{\pi(p)}\partial_r)=\frac{J\partial_r}{r(p)}.$$
Therefore, $\bg(\bnabla_{\partial_r}\partial_r,J\partial_r)=\frac{h_2}{r}dr(J\partial_r)=0$.
Third, for the horizontal lift $Z^*$ of any $Z\in\Gamma(TN)$,
\begin{eqnarray*}
\bg(\bnabla_{\partial_r}\partial_r,Z^*)&=&\frac{1}{2}\Big(\partial_r(\bg(\partial_r,Z^*))+\partial_r(\bg(Z^*,\partial_r))-Z^*(\bg(\partial_r,\partial_r))+\bg([\partial_r,\partial_r],Z^*)\\
&&\phantom{\frac{1}{2}\Big(}-\bg([\partial_r,Z^*],\partial_r)+\bg([Z^*,\partial_r],\partial_r)\Big)\\
&=&\frac{1}{2}\Big(-Z^*(h_2)+2\bg([Z^*,\partial_r],\partial_r)\Big)\\
&=&0.
\end{eqnarray*}
We can already deduce that $\bnabla_{\partial_r}\partial_r=\frac{1}{h_2}\bg(\bnabla_{\partial_r}\partial_r,\partial_r)\partial_r=\frac{h_2'}{2h_2}\partial_r$.
Since $\bnabla J=0$ by construction, we also have $\bnabla_{\partial_r}J\partial_r=J(\bnabla_{\partial_r}\partial_r)=\frac{h_2'}{2h_2}J\partial_r$.
On the other hand, $\bnabla_{J\partial_r}\partial_r=\bnabla_{\partial_r}J\partial_r+[J\partial_r,\partial_r]=\left(\frac{h_2'}{2h_2}+\frac{1}{r}\right)J\partial_r$.\\
Let now $Z^*$ be the horizontal lift of any $Z\in\Gamma(TN)$, then $\bg(\bnabla_{Z^*}\partial_r,\partial_r)=\frac{1}{2}Z^*(\bg(\partial_r,\partial_r))=0$.
Moreover, by $[Z^*,\partial_r]=0$,
\begin{eqnarray*}
\bg(\bnabla_{Z^*}\partial_r,J\partial_r)&=&\frac{1}{2}\Big(Z^*(\bg(\partial_r,J\partial_r))+\partial_r(\bg(J\partial_r,Z^*))-J\partial_r(\bg(Z^*,\partial_r))+\bg([Z^*,\partial_r],J\partial_r)\\
&&\phantom{\frac{1}{2}\Big(}-\bg([\partial_r,J\partial_r],Z^*)+\bg([J\partial_r,Z^*],\partial_r)\Big)\\
&=&\frac{1}{2}\Big(0+0-0-h_2dr(J[Z^*,\partial_r])+\frac{1}{r}\bg(J\partial_r,Z^*)+h_2dr([J\partial_r,Z^*])\Big)\\
&=&\frac{h_2}{2}\cdot(J\partial_r(Z^*(r))-Z^*(J\partial_r(r)))\\
&=&0.
\end{eqnarray*}
For the horizontal lift $Z'^*$ of any $Z'\in\Gamma(TN)$, we have 
\begin{eqnarray*}
\bg(\bnabla_{Z^*}\partial_r,Z'^*)&=&\frac{1}{2}\Big(Z^*(\bg(\partial_r,Z'^*))+\partial_r(\bg(Z'^*,Z^*))-Z'^*(\bg(Z^*,\partial_r))+\bg([Z^*,\partial_r],Z'^*)\\
&&\phantom{\frac{1}{2}\Big(}-\bg([\partial_r,Z'^*],Z^*)+\bg([Z'^*,\partial_r],Z^*)\Big)\\
&=&\frac{1}{2}\partial_r(h_1g_N(Z,Z')\\
&=&\frac{h_1'}{2h_1}\bg(Z^*,Z'^*).
\end{eqnarray*}
Therefore, $\bnabla_{Z^*}\partial_r=\frac{h_1'}{2h_1}Z^*$. 
Note that this remains true for any horizontal vector field since $\bnabla_{\cdot}\partial_r$ is a tensor, it only depends on the pointwise value of vector fields.
As elementary consequences, we have $\bnabla_{Z^*}J\partial_r=J(\bnabla_{Z^*}\partial_r)=\frac{h_1'}{2h_1}JZ^*$ as well as 
$\bnabla_{\partial_r}\overline{Z}=\frac{h_1'}{2h_1}\overline{Z}+[\partial_r,\overline{Z}]$ and $\bnabla_{J\partial_r}\overline{Z}=\frac{h_1'}{2h_1}J\overline{Z}+[J\partial_r,\overline{Z}]$ for any horizontal vector field $\overline{Z}$ on $L$; observe that both brackets $[\partial_r,\overline{Z}]$ and $[J\partial_r,\overline{Z}]$ are horizontal.\\
It remains to compute, for the horizontal lifts $Z^*,Z''^*$ of any vector fields $Z,Z''\in\Gamma(TN)$ and any horizontal vector field $\overline{Z}'$ on $L$,
\begin{eqnarray*}
\bg(\bnabla_{Z^*}\overline{Z}',Z''^*)&=&\frac{1}{2}\Big(Z^*(\bg(\overline{Z}',Z''^*))+\overline{Z}'(\bg(Z''^*,Z^*))-Z''^*(\bg(Z^*,\overline{Z}'))\\
&&\phantom{\frac{1}{2}\Big(}+\bg([Z^*,\overline{Z}'],Z''^*)-\bg([\overline{Z}',Z''^*],Z^*)+\bg([Z''^*,\overline{Z}'],Z^*)\Big)\\
&=&\frac{h_1}{2}\Big(Z^*(g_N(d\pi(\overline{Z}'),Z''))+d\pi(\overline{Z}')(g_N(Z'',Z))-Z''(g_N(Z,d\pi(\overline{Z}')))\\
&&\phantom{\frac{h_1}{2}\Big(}+g_N([Z,d\pi(\overline{Z}')],Z'')-g_N([d\pi(\overline{Z}'),Z''],Z)+g_N([Z'',d\pi(\overline{Z}')],Z)\Big)\\
&=&\overline{g}((\nabla_Z^N d\pi(\overline{Z}'))^*,Z''^*).
\end{eqnarray*}
Here we make the following notation abuse: by $\nabla_Z^N d\pi(\overline{Z}')$ we mean the covariant derivative in $N$ of the vector field $d\pi(\overline{Z}')$ along $Z$ that is obtained by projecting down $\overline{Z}'$ along an integral curve of $Z^*$; this makes sense since, although $d\pi(\overline{Z}')$ is {\sl a priori} not globally defined on $N$, the vector field $\bnabla_{Z^*}Z'^*$ only depends on the value of $\overline{Z}'$ along integral curves of $Z^*$.
On the whole,
\begin{eqnarray*}
\bnabla_{Z^*}\overline{Z}'&=&\frac{1}{h_2}\bg(\bnabla_{Z^*}\overline{Z}',\partial_r)\partial_r+\frac{1}{h_2}\bg(\bnabla_{Z^*}\overline{Z}',J\partial_r)J\partial_r+(\nabla_Z^N d\pi(\overline{Z}'))^*\\
&=&-\frac{1}{h_2}\bg(\overline{Z}',\bnabla_{Z^*}\partial_r)\partial_r-\frac{1}{h_2}\bg(\overline{Z}',\bnabla_{Z^*}J\partial_r)J\partial_r+(\nabla_Z^N d\pi(\overline{Z}'))^*\\
&=&-\frac{h_1'}{2h_1h_2}\left(\bg(Z^*,\overline{Z}')\partial_r+\bg(JZ^*,\overline{Z}')J\partial_r\right)+(\nabla_Z^N d\pi(\overline{Z}'))^*.
\end{eqnarray*}
In particular the computations show that $J\partial_r$ is a Killing vector field on $(L,\bg)$.
We are now ready to prove Proposition \ref{l:HessianeqonCalabi}.
Let $u$ be as in the assumptions.
Then $du=\partial_ru\cdot dr=u'(r)\cdot dr$ and therefore $\bnabla_Xdu=X(u'(r))dr+u'(r)\bnabla_Xdr$ for any $X\in TL$.
Using the splitting $X=\frac{1}{h_2}\bg(X,\partial_r)\partial_r+\frac{1}{h_2}\bg(X,J\partial_r)J\partial_r+X^H$, where $X^H$ is the pointwise orthogonal projection of $X$ onto the horizontal subbundle $\mathcal{H}L$, we obtain
\begin{eqnarray*}
\bnabla_X^2u&=&X(u'(r))\frac{1}{h_2}\partial_r+u'(r)\bnabla_X(\frac{1}{h_2}\partial_r)\\
&=&\frac{1}{h_2}\left(X(u'(r))-\frac{X(h_2)}{h_2}u'(r)\right)\partial_r+\frac{u'(r)}{h_2}\Big(\frac{1}{h_2}\bg(X,\partial_r)\bnabla_{\partial_r}\partial_r+\frac{1}{h_2}\bg(X,J\partial_r)\bnabla_{J\partial_r}\partial_r+\bnabla_{X^H}\partial_r\Big)\\
&=&\frac{1}{h_2}\left(X(u'(r))-\frac{X(h_2)}{h_2}u'(r)\right)\partial_r\\
&&+\frac{u'(r)}{h_2}\Big(\frac{h_2'}{2h_2^2}\bg(X,\partial_r)\partial_r+\frac{1}{h_2}(\frac{h_2'}{2h_2}+\frac{1}{r})\bg(X,J\partial_r)J\partial_r+\frac{h_1'}{2h_1}X^H\Big)\\
&=&\frac{1}{h_2}\Big(X(u'(r))+\frac{u'(r)}{h_2}(\frac{h_2'}{2h_2}\bg(X,\partial_r)-X(h_2))\Big)\partial_r+\frac{u'(r)}{h_2^2}\left(\frac{h_2'}{2h_2}+\frac{1}{r}\right)\bg(X,J\partial_r)J\partial_r\\
&&+\frac{u'(r)h_1'}{2h_1h_2}X^H\\
&=&\frac{1}{h_2^2}\Big(u''(r)-\frac{u'(r)h_2'}{2h_2}\Big)\bg(X,\partial_r)\partial_r+\frac{u'(r)}{h_2^2}(\frac{h_2'}{2h_2}+\frac{1}{r})\bg(X,J\partial_r)J\partial_r+\frac{u'(r)h_1'}{2h_1h_2}X^H.
\end{eqnarray*}
As a consequence, for the horizontal lift $X^*$ of any $X\in\Gamma(TN)$,
\begin{equation}\label{eq:hessianucalabi}
\left\{\begin{array}{ll}\bnabla_{\partial_r}^2u&=\frac{1}{h_2}\Big(u''-\frac{u'h_2'}{2h_2}\Big)\partial_r\\
&\\
\bnabla_{J\partial_r}^2u&=\frac{u'}{h_2}(\frac{h_2'}{2h_2}+\frac{1}{r})J\partial_r\\
&\\
\bnabla_{X^*}^2u&=\frac{u'h_1'}{2h_1h_2}X^*\end{array}\right.
\end{equation}
All three vectors $\partial_r,J\partial_r,X^*$ are already pointwise eigenvectors of $u$ and the only condition to be satisfied is the identity $u''-\frac{u'h_2'}{2h_2}=(\frac{h_2'}{2h_2}+\frac{1}{r})u'$, that is,
$$u''-(\frac{h_2'}{h_2}+\frac{1}{r})u'=0.$$
This is a first-order linear ODE in $u'$ whose solution is $u'(r)=C\cdot h_2(r)\cdot r$ for some real constant $C$, i.e. $u(r)=u(0)+C\cdot\int_0^rs\cdot h_2(s)ds$.
Mind that we make here a light notation abuse since we consider the function $r$ {\sl a priori} away from $0$.
If $C=0$, then $u$ is constant which is the trivial case.
Otherwise, up to replacing $u$ by $\frac{u-u(0)}{C}$ (which does not modify the assumption on the Hessian on $u$ nor on its eigenspaces), we may assume that $u(0)=0$ and that $C=1$, so that
$$u(r)=\int_0^rs\cdot h_2(s)ds.$$
Together with (\ref{eq:hessianucalabi}) we obtain the expressions for $\lambda$ and $\mu$.
Note that, coming back to the original notations for $h_1,h_2$, which are defined as functions of $r^2$, we obtain via the change of variable $\tau:=s^2$ the formula
$$u(r)=\frac{1}{2}\int_0^{r^2}h_2(\tau)d\tau.$$
Furthermore, in case $l\neq0$, we obtain the alternative expression $u(r)=-\frac{1}{l}h_1(r^2)+\mathrm{cst}$ for $u$.
\findemo

It must be noticed that the assumptions on $u$ are strong but that they will turn out to be fulfilled in one of the two main cases of our Obata-like equation.
We also point out that the conclusions of Lemma \ref{l:HessianeqonCalabi} might hold without the assumption $d^Hu=0$ in the sense that that condition might be implied by the other assumptions.

\subsection{Proof of Theorem \ref{t:cardnormalbundle}}\label{ss:proofmaintheorem}

We come to the proof of Theorem \ref{t:cardnormalbundle}.
Note first that the assumption of $\wit{J}$-invariance of $\wnabla^2u$, which is equivalent to the holomorphicity of $\wnabla u$ (or of $\wit{J}(\wnabla u)$), ensures the critical set $N:=\{x\in \wit{M}\,/\,(\wnabla u)_x=0\}$ to have empty interior since $\wit{M}$ is connected and $u$ is assumed to be nonconstant: any analytic vector field vanishing on a nonempty open subset must vanish identically on any connected component meeting that open set.
Since $\wnabla^2u$ is furthermore symmetric, the holomorphic vector field $\wit{J}(\wnabla u)$ is also Killing and therefore $N$ is a totally geodesic K\"ahler submanifold of dimension $2k<2n$ of $\wit{M}$, see e.g. \cite[Sec. II.5]{Kocomplexgeom}.
Because $u$ is furthermore assumed to be convex, the submanifold $N$ must be the set of absolute minima of $u$ and must be connected: considering namely any geodesic joining any two critical points of $u$, the function $u$ along that geodesic is convex with two critical points, therefore the corresponding critical value is a minimum and the function is constant between the two critical points by convexity.
Moreover, because $\wnabla u$ is analytic, so is $\wnabla^2 u$, so that -- at least locally around every point -- there is an analytic family of eigenvalues and eigenspaces associated to $\wnabla^2 u$, in particular $\lambda$ and $\mu$ as well as their eigenspaces extend through $N$.
Since $TN=\ker\left(\wnabla(\wit{J}\wnabla u)\right)=\ker\left(\wnabla^2u\right)$ (see e.g. \cite[Sec. II.5]{Kocomplexgeom}), both $\lambda$ and $\mu$ cannot simultaneously vanish at any point of $N$, otherwise $u$ would be constant.
In the same way, the vanishing set of $u$ necessarily has empty interior. 
We fix a regular value $u_0$ of $u$ and set $M:=u^{-1}(\{u_0\})$.
We denote by $\nu:=\frac{\wnabla u}{|\wnabla u|}$ the normalized gradient vector field of $u$ on its regular set and by $\xi:=-\wit{J}\nu$.
We let $\hat{\xi}:=\alpha\xi$ and $\hat{g}:=\alpha^{-2}g_{\hat{\xi}}\oplus\beta(x)^{-2}g_{\hat{\xi}^\perp}$ for some positive $\alpha\in\R$ and basic function $\beta$ on $M$ which is constant when $n>2$; both $\alpha$ and $\beta$ will be fixed later on.
The flow $F$ of $\nu$ is defined at least on $(-\varepsilon,\varepsilon)\times M$ for some $\varepsilon>0$ since $|\wnabla u|$ is constant along any regular level hypersurface of $u$; in other words, we get a neighbourhood $U$ of $M$ in $\wit{M}\setminus N$ that is uniformly thick around $M$.
We write $U$ as a doubly-warped product as in Theorem \ref{t:cardwp2eigenv}.
First, if $\mu$ vanishes at one (regular) point of $u$ in $U$, then $\mu$ vanishes identically on $U$ by Theorem \ref{t:cardwp2eigenv} and hence on $\wit{M}$ by analyticity.
If $\mu=0$ on $\wit{M}$, then in particular $\mu_{|_N}=0$ .
But if $\mu_{|_N}=0$, then because the eigenvalue $\lambda$ cannot vanish at any point of $N$, necessarily $\dim(N)=2n-2$ and $T^\perp N=\ker(\wnabla^2u-\lambda\cdot\mathrm{Id})$ hold.
In case $\mu$ does not vanish identically on $\wit{M}$, then again by the arguments brought up in the proof of Theorem \ref{t:cardwp2eigenv}, we have, up to suitably choosing $\alpha,\beta$ and replacing $u$ by $au+b$ for $a,b\in\R$ with $a\neq0$, that $u>0$ and $\mu=\frac{|\wnabla u|^2}{2u}$ on $U$.
Note that $a$ must be positive by the convexity assumption on $u$.
Also, as in the proof of Theorem \ref{t:cardwp2eigenv}, it can be shown that the (closed nonempty) set of points of $\wit{M}\setminus N$ where
\begin{equation}\label{eq:exprmu}2u\mu=|\wnabla u|^2\end{equation}
is satisfied is open in $\wit{M}\setminus N$.
Because the critical submanifold $N$ has codimension at least $2$, the open set $\wit{M}\setminus N$ is connected and dense in $\wit{M}$, therefore (\ref{eq:exprmu}) is fulfilled on $\wit{M}\setminus N$ and thus on $\wit{M}$.

If $u$ vanishes nowhere, in which case $\mu=\frac{|\wnabla u|^2}{2u}$ holds on $\wit{M}$, then in particular $\mu_{|_N}=0$ and, by the above, $N$ is a $(2n-2)$-dimensional totally geodesic K\"ahler submanifold  of $\wit{M}$ with normal bundle $\ker(\wnabla^2u-\lambda\cdot\mathrm{Id})$.
If now $u(x)=0$ for some $x\in\wit{M}$, then $(\wnabla u)_x=0$ must hold, in particular $x\in N$.
We next show that $\lambda(x)=\mu(x)$, in particular $(\wnabla^2u)_x=\lambda(x)\cdot\mathrm{Id}\neq0$.
Namely let $v\in T_x\wit{M}$ be any nonzero vector and $\gamma_v\colon \R\to\wit{M}$ be any smooth curve with $\gamma_v(0)=x$ and $\dot{\gamma}_v(0)=v$.
Let $y_v:=u\circ\gamma_v\colon \R\to \R$, then $y_v$ is smooth with $y_v(0)=0$ as well as $y_v'=\langle\wnabla u,\dot{\gamma}_v\rangle$.
Equation (\ref{eq:exprmu}) translates as $(\mu\circ\gamma_v)\cdot y_v= \frac{1}{2}|\wnabla u|^2\circ\gamma_v$.
But
$$\frac{1}{2}\left(|\wnabla u|^2\circ\gamma_v\right)'=\langle\wnabla_{\dot{\gamma}_v}\wnabla u,\wnabla u\rangle=\langle\wnabla_{\wnabla u}^2u,\dot{\gamma}_v\rangle=(\lambda\circ\gamma_v)\langle\wnabla u,\dot{\gamma}_v\rangle=(\lambda\circ\gamma_v)y_v',$$
so that 
$$(\mu\circ\gamma_v)'\cdot y_v+(\mu\circ\gamma_v)\cdot y_v'=\left((\mu\circ\gamma_v)\cdot y_v\right)'=\frac{1}{2}\left(|\wnabla u|^2\circ\gamma_v\right)'=(\lambda\circ\gamma_v)y_v'.$$
In other words, $y_v$ satisfies the \emph{linear} first-order ODE
$$(\mu\circ\gamma_v-\lambda\circ\gamma_v)\cdot y_v'+(\mu\circ\gamma_v)'\cdot y_v=0.$$
By contradiction, if $\lambda(x)\neq\mu(x)$, then $\mu\circ\gamma_v-\lambda\circ\gamma_v$ would not vanish on an open interval $I_v$ about $0$ and therefore $y_v'=-\frac{(\mu\circ\gamma_v)'}{\mu\circ\gamma_v-\lambda\circ\gamma_v}\cdot y_v$ would hold on $I_v$, which with $y_v(0)=0$ would yield $y_v=0$ on $\R$.
But then $0=y_v''(0)=\langle\wnabla_v^2u,v\rangle$ would hold and for any nonzero $v\in T_x\wit{M}$.
In turn this would imply $(\wnabla^2u)_x=0$ and in particular $\lambda(x)=\mu(x)=0$, contradiction.
Therefore $\lambda(x)=\mu(x)$ must hold and therefore $(\wnabla^2u)_x=\lambda(x)\cdot\mathrm{Id}$; the eigenvalue $\lambda(x)=\mu(x)$ then cannot vanish as we mentioned at the beginning of the proof.\\
To sum up, we have shown that, in the case where $\mu\neq0$ on $\wit{M}$, for any critical point $x$ of $u$, either $u(x)\neq0$ and then $\mu(x)=0$, $\lambda(x)\neq0$, $\dim(N)=2n-2$ and $T^\perp N=\ker(\wnabla^2u-\lambda\cdot\mathrm{Id})$; or $u(x)=0$ and then $\mu(x)=\lambda(x)\neq0$ and hence $N=\{x\}$.
Note that we implicitly use the fact that $N$ is connected, which is true if $u$ is convex, otherwise the preceding argument would apply to each connected component of $u$.
Note also that in any case $\lambda$ cannot vanish at any critical point of $u$.\\

Next we show that $u$ only depends on the geodesic distance to $N$, which in particular relies on the fact that $N$ is connected.
Let $\gamma\colon\R\to\wit{M}$ be any geodesic with $\gamma(0)\in N$ and $\dot{\gamma}(0)\in T_{\gamma(0)}^\perp N$.
Note that, in all cases, $\dot{\gamma}(0)\in\ker(\wnabla^2u-\lambda\cdot\mathrm{Id})$ holds.
We start with the case where $\dim_\R(N)=2n-2$.
Note that $\ker(\wnabla^2u-\lambda\cdot\mathrm{Id})$ is a totally geodesic integrable distribution in $\wit{M}$: this is a straightforward consequence of the equation $\wnabla^2u=\lambda\cdot\mathrm{Id}_{\mathrm{Span}(\wnabla u,J\wnabla u)}+\mu\cdot\mathrm{Id}_{\mathrm{Span}(\wnabla u,J\wnabla u)^\perp}$ since, outside $N$, both $\wnabla u$ and $J\wnabla u$ are $\lambda$-eigenvectors of $\wnabla^2u$.
In particular, if $\dot{\gamma}(0)\in\ker(\wnabla^2u-\lambda\cdot\mathrm{Id})$, then $\gamma$ must remain tangent to $\ker(\wnabla^2u-\lambda\cdot\mathrm{Id})$, that is, there exist coefficients $a(t),b(t)$ with $\dot{\gamma}(t)=a(t)\nu_{\gamma(t)}+b(t)\xi_{\gamma(t)}$ for all $t\neq0$.
Now the function $\langle\dot{\gamma},\xi\rangle$ is equal to $-\frac{1}{|(\wnabla u)\circ \gamma|}\langle\dot{\gamma},J\wnabla u\rangle$ outside $0$ and $\langle\dot{\gamma},J\wnabla u\rangle$ is constant on $\R$: its derivative is given by $\langle\frac{\wnabla\dot{\gamma}}{dt},J\wnabla u\rangle+\langle\dot{\gamma},\wnabla_{\dot{\gamma}}J\wnabla u\rangle=0$ because of $J(\wnabla u)$ being Killing.
Since $\langle\dot{\gamma}(0),J\wnabla u\rangle=0$ (by $(\wnabla u)_{\gamma(0)}=0$), we deduce that $\langle\dot{\gamma},J\wnabla u\rangle=0$ and hence $\langle\dot{\gamma},\xi\rangle=0$ on $\R\setminus\{0\}$.
This shows $b(t)=0$ and hence $\dot{\gamma}(t)=a(t)\nu_{\gamma(t)}$ for all $t\neq0$.
Therefore, we have proved that any geodesic leaving $N$ normally must become immediately tangent to $\nu$.\\
{In the case where $N=\{x\}$, we show the same fact but write a different proof that is strongly inspired by the proof of \cite[Lemma 17]{Kuehnel}.
We fix a convex open neighbourhood $U$ of $x$, in particular there is a unique geodesic joining any point $y\in U$ with $x$.
Since $\wnabla^2u$ is positive definite at $x$ and hence in a neighbourhood of $x$, any level hypersurface associated to a value of $u$ closed to its minimum value must be compact -- it is diffeomorphic to the sphere $\mathbb{S}^{2n-1}$ by the Morse lemma -- and therefore all level hypersurfaces are compact since they are all diffeomorphic.
Let $\gamma\colon\R\to \wit{M}$ be any geodesic with $\gamma(0)=x$.
Let $t_0>0$ be so small that $\gamma(t_0)\in U$ and even such that the smooth level hypersurface $M_{u(\gamma(t_0))}:=u^{-1}(\{u(\gamma(t_0))\})$ lies entirely inside $U$.
Recall that, by assumption, both distributions $\ker(\wnabla^2u-\lambda\mathrm{Id})$ and $\ker(\wnabla^2u-\mu\mathrm{Id})$ are analytic and defined on the whole $\wit{M}$, even at the critical point $x$ of $u$.
Moreover, the $2$-dimensional distribution $\ker(\wnabla^2u-\lambda\mathrm{Id})$ is integrable to a totally geodesic two-dimensional foliation.
In particular, there is a unique such leaf $S$ through $x$.
Let now $\gamma_1\colon\R\to\wit{M}$ be a geodesic with $\gamma_1(0)=x$ and $\dot{\gamma_1}(0)\in T_xS$, in particular $\gamma_1(t)\in S$ for very $t\in\R$.
Note that, since $\dot{\gamma}_1(t)$ is at every $t\neq0$ a linear combination of $\xi$ and $\nu$, the same argument as above shows that actually $\dot{\gamma}_1(t)$ must be tangent to $\nu$ for every $t\neq0$.
We can assume w.l.o.g. that $\dot{\gamma}_1(t)=\nu$ for all $t\neq0$.
Let $s_0\in\R$ be such that $\gamma_1(s_0)\in M_{u(\gamma(t_0))}$ (there is a unique such $s_0$ since the convex function $u$ must be monotonously increasing and unbounded along the geodesic $\gamma_1$ from $x$).
Let $\gamma_2\colon\R\to\wit{M}$ be the $\nu$-geodesic with $\gamma_2(s_0)=\gamma(t_0)$.
Then as in the proof of \cite[Lemma 17]{Kuehnel} and using the fact that $\gamma_1(t),\gamma_2(t)$ lie on the same level hypersurface of $u$ for all $t>0$, we write, for every $t>0$,
\[d_{\wit{M}}(\gamma_1(t),\gamma_2(t))\leq d_{M_{u(\gamma_2(t))}}(\gamma_1(t),\gamma_2(t))=(\rho\rho')(t)d_{u(\gamma(t_0))}(\gamma_1(s_0),\gamma_2(s_0)),\]
so that $d_{\wit{M}}(\gamma_1(t),\gamma_2(t))\buil{\longrightarrow}{t\to0}0$ because of $\rho'(t)\buil{\longrightarrow}{t\to0}0$.
This shows that $\gamma_1(0)=\gamma_2(0)$, in particular $\gamma_2$ is another geodesic joining $x$ and $\gamma(t_0)$ and thus $\gamma_2=\gamma_1=\gamma$ because of both $\gamma,\gamma_2$ lying in $U$ (use the fact that $u$ is proper, at least in $U$).
This proves that $\gamma$ gets immediatly tangent to $\nu$ when leaving $x$.\\

The fact that every geodesic leaving $N$ normally (any geodesic when $N=\{x\}$) gets immediatly tangent to $\nu$ first implies that any integral curve of $\nu$ meets $N$ normally.
For given any $x\in\wit{M}\setminus N$, let $\gamma$ be the $\nu$-geodesic starting at $x$.
Let $\gamma_{\min}$ be a minimizing geodesic from $x$ to the closest point to $x$ on $N$ (that point exists since $N$ is closed, the distance function is continuous and $\wit{M}$ is locally compact), then $\gamma_{\min}$ meets $N$ orthogonally and therefore is tangent to $\nu$, so that it must coincide with $\gamma$ up to reparametrization and in particular $\gamma$ must meet $N$ normally and at a closest point to $x$.\\
Let now $x,y\in\wit{M}$ be any two points with $d(x,N)=d(y,N)$.
Let $\gamma_x,\gamma_y$ be the $\nu$-geodesic starting from $x$ and $y$ respectively.
As we have just shown, $\gamma_x$ (resp. $\gamma_y$) meets $N$ normally at a closest point $x'$ (resp. $y'$), that is, $x',y'\in N$ with $d(x,x')=d(x,N)$ (resp. $d(y,y')=d(y,N)$).
But since $\gamma_x$ and $\gamma_y$ have unit speed, they are parametrized by arc-length, so that $d(x,x')$ (resp. $d(y,y')$) coincides with the time $\gamma_x$ (resp. $\gamma_y$) takes from $x'$ to $x$ (resp. from $y'$ to $y$).
Because $d(x,x')=d(y,y')$, both times are equal.
But since $N$ is connected, then $u_{|_N}$ is constant; { more precisely, we claim that $u\circ\gamma_z$ only depends on time and not on the starting point $z\in N$.
Namely recall that this property holds true if we start from a \emph{regular} hypersurface of $u$, as can be deduced from Proposition \ref{p:warpedprodu}.
Let $z,z'\in N$ be any two points and consider normal geodesic $\gamma_z,\gamma_{z'}$ starting at $z,z'$ respectively.
Fix some $t_0\in\R$ and consider the regular real hypersurface $M_0:=u^{-1}(\{u\circ\gamma_z(t_0)\})\subset\wit{M}$.
The geodesic $\gamma_{z'}$ must meet $M_0$ because $u$ is monotonously increasing and unbounded along any normal geodesic (recall that $u$ is assumed to be convex and nonconstant). 
Let $t_0'\in\R$ be such that $\gamma_{z'}(t_0')\in M_0$.
We claim that $t_0=t_0'$: for following $\gamma_{z'}$ and $\gamma_z$ backward from $\gamma_{z'}(t_0'),\gamma_z(t_0)$ respectively we have $u\circ\gamma_{z'}(t_0'-s)=u\circ\gamma_z(t_0-s)$ for all $s$ by Proposition \ref{p:warpedprodu}, in particular $u$ reaches its value $u_{|_N}=\min(u)$ after the time $-t_0$, therefore $u(\gamma_{z'}(t_0'-t_0))=\min(u)$ and hence $t_0'-t_0=0$.
On the whole, we obtain that, if $d(x,N)=d(y,N)$, then $u(x)=u(y)$.


Therefore, $u$ only depends on the distance to $N$, in other words the level hypersurfaces of $u$ coincide with the distance tubes around $N$.\\

Next we show that, because of $u$ being convex by assumption, the normal exponential map $\exp^\perp\colon T^\perp N\to\wit{M}$ is a global diffeomorphism.
The existence of such a diffeomorphism actually follows from \cite[Theorem B]{GreeneShiohamadiff81} using only convexity of $u$, see also \cite{GreeneShiohamatop81} for the topological implications of convexity.\\
As we have already seen above, each geodesic $\gamma$ of $\wit{M}$ with $\dot{\gamma}(0)\in T^\perp N$ must be a ray, that is, it must minimize distance to $\gamma(0)$ for all times.
Moreover, $\gamma$ cannot meet $N$ twice since $u$ is strictly monotonous along $\gamma$.
This fact ensures that $\exp^\perp\colon T^\perp N\to\wit{M}$ is a diffeomorphism onto its range, which is hence an open set of $\wit{M}$.
The sujectivity of $\exp^\perp$ is straightforward consequence of all $\nu$-geodesics meeting $N$ normally.
On the whole, the map $\exp^\perp\colon T^\perp N\to \wit{M}$ is a global diffeomorphism.\\

Next we describe the metric structure of $\wit{M}$, that is, we pull $\tilde{g}$ back onto $T^\perp N$ via $\exp^\perp$.
First, because every normal geodesic leaving $N$ becomes immediately tangent to $\nu$, the Gau\ss{} lemma implies that the normal exponential  $\exp^\perp$ sends the usual normalized position vector field in each fibre onto $\nu$, see e.g. \cite[Cor. 2.14]{Graytubes}.
This has the following fundamental consequence that the normal exponential map identifies with the flow of $\nu$:
denoting namely by $\mathrm{proj}\colon\wit{M}\to N$ the map sending some $x$ onto the unique intersection point of a $\nu$-geodesic through $x$ with $N$ in case $x\in\wit{M}\setminus N$ and $\mathrm{proj}(x)=x$ if $x\in N$, we have $\mathrm{proj}\circ\exp^\perp=\pi$, where $\pi\colon T^\perp N\to N$ is the standard footpoint map of the normal bundle $T^\perp N$ of $N$.

Assume now that $\dim_\R(N)=2n-2$, case that happens as soon as $\mu$ vanishes somewhere.
Consider $(\exp^\perp)^*\wit{J}$ instead of the original almost-Hermitian structure on the total space of $T^\perp N\to N$.
Note that $(\exp^\perp)^*\wit{J}$ coincides with the restriction of the original complex structure on the horizontal distribution $\ker(\wnabla^2 u-\mu\mathrm{Id})$ since the flow of $\nu$ preserves $\wit{J}$ along that distribution as we have seen in the proof of Theorem \ref{t:cardwp2eigenv}.
However, $(\exp^\perp)^*\wit{J}$ does not coincide with $\wit{J}$ along $\ker(\wnabla^2u-\lambda\mathrm{Id})$ since $\mathcal{L}_\nu\wit{J}$ does not vanish along that distribution.
As for $\mathrm{Span}(\wnabla u,J\wnabla u)^\perp=\ker(\wnabla^2 u-\mu\cdot\mathrm{Id})$, we have already shown that, at least on the open subset of regular points of $u$, this distribution is preserved by the flow of $\nu$; moreover, the induced metric from $\wit{g}$ is preserved up to multiplication by a factor depending on $t$ only, i.e. on the value of $u$.
Therefore, taking the limit when getting close to $N$, we obtain that the map $\mathrm{proj}$ identifies $\wit{g}_{|_{\ker(\wnabla^2 u-\lambda\mathrm{Id})}}$ with $h_1\cdot\pi^*g_N$ for some positive function $h_1$ depending only on the distance function to $N$ and where $g_N$ is the restriction of $\wit{g}$ onto $N$.
On the whole, the manifold $(T^\perp N,(\exp^\perp)^*\wit{g},(\exp^\perp)^*\wit{J})$ is K\"ahler with metric given by 
\[(\exp^\perp)^*\wit{g}=h_1\cdot\pi^* g_N+g_V,\]
where $g_V:=(\exp^\perp)^*\wit{g}_{|_{T^\perp N}}$ is the vertical part of the metric.
Note that in every fibre $g_V$ is rotationnally invariant.
Since $T^\perp N$ has rank $2$, the metric $g_V$ is conformal on every fibre to the flat Euclidean metric.
Because both $g_V$ and that flat metric are rotationnally invariant, so is the conformal factor, therefore it is a function of the fibre length.
The complex structure $(\exp^\perp)^*\wit{J}$ is orthogonal w.r.t. $(\exp^\perp)^*\wit{g}$ and therefore also w.r.t. to the conformal flat metric.
But on flat $\R^2$ there is only one complex structure up to sign, therefore $(\exp^\perp)^*\wit{J}$ coincides with the standard complex structure $J$ on flat $\R^2$ up to sign.
Now that sign must be equal to $1$ for every fibre since along the zero section the map $\exp^\perp$ identifies both the complex structure of the metric and $\wit{J}$.
Defining $r$ to be the norm w.r.t. the flat metric in each fibre, the normal exponential map defines a biholomorphic isometry from $(T^\perp N,(\exp^\perp)^*\wit{g},(\exp^\perp)^*\wit{J})$ to $(\wit{M},\wit{g},\wit{J})$, where
\[(\exp^\perp)^*\wit{g}=h_1(r)\cdot\pi^* g_N+h_2(r)\cdot(dr\otimes dr+d^c r\otimes d^c r)\]
and $d^c r:=-dr\circ J$.
But this is exactly the form we want for the Calabi construction.
It follows from Proposition \ref{l:HessianeqonCalabi}, whose assumptions are thus all fulfilled, that there exists a constant $l$ such that $h_1'(r)=-lrh_2(r)$, and that $u(r)=\int_0^{r}s\cdot h_2(s)ds$ up to replacing $u$ by $bu+b'$ with $b,b'\in\R$.
Note that necessarily $b>0$ because of the convexity of $u$.
It also directly shows that, in case $l=0$, which corresponds to the case where $\mu=0$ by Proposition \ref{l:HessianeqonCalabi}, the normal bundle $T^\perp N\to N$ is flat and that in case $l\neq0$ the K\"ahler manifold $N$ must be Hodge.
If $l=0$, the manifold $\wit{M}$ must actually be locally biholomorphically isometric to the product $(\R^2\times N,h_2(r)\cdot(dr\otimes dr+ d^c r\otimes d^c r)\oplus g_N)$; that identification is global if e.g. $\wit{M}$ is assumed to be simply-connected.

Now we look at the metric in the case where $N=\{x\}$ and hence $\lambda(x)=\mu(x)\neq0$.
By convexity of $u$, we necessarily have $\lambda(x)=\mu(x)>0$.
Note also that, because $N$ is reduced to a point, the distance tubes around $N$ are diffeomorphic to spheres, that is, $M=\mathbb{S}^{2n-1}$ with possibly nonround metric.
Thus, outside $x$, the flow of $\nu$ identifies $\wit{M}\setminus\{x\}$ with some doubly warped product manifold $(0,\infty)\times\mathbb{S}^{2n-1}$ with metric $dt^2+\rho(t)^2(\rho'(t)^2\hat{g}_{\hat{\xi}}+\hat{g}_{\hat{\xi}}^\perp)$, where $(\mathbb{S}^{2n-1},\hat{g},\hat{\xi})$ is Sasaki.
Up to translating $u$, we have seen that we may assume that $\mu=\frac{|\wnabla u|^2}{2u}$; up to multiplying $u$ by a constant, we may also assume that $\lambda(x)=\mu(x)=2$.
We compute, using Gau\ss{} formula, the sectional curvatures of $(M,g_t=\rho(t)^2(\rho'(t)^2\hat{g}_{\hat{\xi}}+\hat{g}_{\hat{\xi}}^\perp))$.
The second fundamental form of $(M,g_t)$ in $\wit{M}$ is given for all $X,Y\in\Gamma(TM)$ by 
$$II(X,Y)=-\wit{g}(\wnabla_X\nu,Y)\nu=-\wit{g}(\wnabla_X(\frac{\wnabla u}{|\wnabla u|}),Y)\nu=-\frac{1}{|\wnabla u|}\wit{g}(\wnabla_X^2u,Y)\nu,$$
hence $II(\xi,\xi)=-\frac{\lambda}{|\wnabla u|}\nu$, $II(\xi,Z)=0$ and $II(Z,Z')=-\frac{\mu}{|\wnabla u|}\wit{g}(Z,Z')\nu$ for all $Z,Z'\in\Gamma(Q)$.
The formula relating the sectional curvatures is given for all linearly independent vectors $X,Y\in TM$ by
$$K^M(X,Y)=K^{\wit{M}}(X,Y)+\frac{\wit{g}(II(X,X),II(Y,Y))-|II(X,Y)|^2}{|X^2|\cdot|Y^2|-\wit{g}(X,Y)^2}.$$
Therefore, 
$K^M(Z,Z')=K^{\wit{M}}(Z,Z')+\frac{\mu^2}{|\wnabla u|^2}$.
Now we relate the sectional curvatures of $(M,g_t)$ and of $(M,\hat{g}=\hat{g}_{\hat{\xi}}+\hat{g}_{\hat{\xi}}^\perp)$.
We have 
$$K^M(Z,Z')=\frac{1}{\rho(t)^2}\left(K^{\hat{M}}(Z,Z')-3(\rho'(t)^2-1)\hat{g}(JZ,Z')^2\right),$$
where here $\{Z,Z'\}$ is chosen to be $\hat{g}$-orthonormal.
Note that, since $(M,\hat{g},\hat{\xi})$ is Sasaki, the sectional curvature of $(M,\hat{g},\hat{\xi})$ in the direction of $\hat{\xi}$ and $Z$ is $1$.
We can deduce that
\begin{eqnarray*}
K^{\wit{M}}(Z,Z')&=&K^M(Z,Z')-\frac{\mu^2}{|\wnabla u|^2}\\
&=&\frac{1}{\rho(t)^2}\left(K^{\hat{M}}(Z,Z')-3(\rho'(t)^2-1)\hat{g}(JZ,Z')^2\right)-\frac{|\wnabla u|^2}{4u^2}\textrm{ with }\frac{|\wnabla u|^2}{4u^2}=(\frac{\rho'}{\rho})^2\\
&=&\frac{1}{\rho(t)^2}\left(K^{\hat{M}}(Z,Z')-3(\rho'(t)^2-1)\hat{g}(JZ,Z')^2-\rho'(t)^2\right).
\end{eqnarray*}
As in the proof of \cite[Lemma 18]{Kuehnel}, we argue as follows: letting $t\to0$ -- i.e, getting close to $N$ along any $\nu$-geodesic -- the l.h.s. of the above identity for sectional curvature has a well-defined finite limit since $\wit{g}$ is smooth and each plane considered has a limit plane along $N$.
This forces the r.h.s. also to have a finite limit.
But $\rho(t)\buil{\longrightarrow}{t\to0}0$, so that necessarily 
$K^{\hat{M}}(Z,Z')-3(\rho'(t)^2-1)\hat{g}(JZ,Z')^2-\rho'(t)^2\buil{\longrightarrow}{t\to0}0$, i.e. $K^{\hat{M}}(Z,Z')\buil{\longrightarrow}{t\to0}3(\rho'(0)^2-1)\hat{g}(JZ,Z')^2+\rho'(0)^2$.
Since $\hat{g}$ does not depend on $t$, this means that $K^{\hat{M}}(Z,Z')=3(\rho'(0)^2-1)\hat{g}(JZ,Z')^2+\rho'(0)^2$.
With $2=\lambda(x)=(\rho^2)''(0)=2\rho'(0)^2+2\rho(0)\rho''(0)=2\rho'(0)^2$, we have $\rho'(0)^2=1$ and hence $K^{\hat{M}}(Z,Z')=1$.
Since on the other hand $K^{\wit{M}}(\hat{\xi},Z)=1$ holds because of $(\mathbb{S}^{2n-1},\hat{g},\hat{\xi})$ being Sasaki, we can conclude that $(\mathbb{S}^{2n-1},\hat{g})$ has constant sectional curvature $1$.
The oddness of $\rho$ is imposed by the smoothness of the metric, see for instance proof of \cite[Lemma 18]{Kuehnel}.
This ends the proof of Theorem \ref{t:cardnormalbundle}.


\begin{erem}\label{r:Berger}
{\rm Note that, if we do not rescale $u$, then we end up with a Berger metric on $\mathbb{S}^{2n-1}$, in particular we obtain again the round metric by further suitable rescaling along each factor.
}
\end{erem}

The question remains open how Theorem \ref{t:cardnormalbundle} can be generalized by deleting the convexity assumption on $u$, i.e. only assuming $\wnabla^2u$ to be everywhere $J$-invariant and to have pointwise two eigenvalues, one of which with $\wnabla u$ as an eigenvector.
The idea would be to recover the results of \cite{MolzonPinneyMortensen93,RanjSant97,Santhanam07} in a unified framework.
In that case, the normal exponential along the critical set of $u$, which has no more reason to be connected, cannot be used any longer as an identification.
This is the object of future work.\\

{\bf Acknowledgment:} This long-term project benefited from the generous support of the Universities of Stuttgart, Lorraine, Regensburg -- in particular from the Johannes-Kepler-Zentrum f\"ur Mathematik -- and from the conference \emph{Riemann and K\"ahler geometry} held at IMAR (Bucharest) between April 15-19, 2019.
We would like to thank especially Sergiu Moroianu and the GDRI Eco-Math for his support and interest in our project.
The second named author would also like to thank the Humboldt Foundation as well as the German Academic Exchange Service (DAAD) for their support.
{\sl Last but not the least}, we are very grateful to Paul-Andi Nagy for pointing out to us the relation of our work to other classification results on K\"ahler manifolds.

\providecommand{\bysame}{\leavevmode\hbox to3em{\hrulefill}\thinspace}

\end{document}